\documentclass[12pt]{article}

\usepackage[english]{babel}
\usepackage{amssymb,a4wide}
\usepackage{amsmath,xypic,epsf,mathrsfs,theorem}
\usepackage[all]{xy}

\newtheorem{theorem}{Theorem}[subsection]
\newtheorem{lemma}[theorem]{Lemma}
\newtheorem{corollary}[theorem]{Corollary}

\newtheorem{proposition}[theorem]{Proposition}

{\theorembodyfont{\upshape}
\newtheorem{definition}[theorem]{Definition}
\newtheorem{remark}[theorem]{Remark}
\newtheorem{example}[theorem]{Example}
}

\numberwithin{equation}{section}
\numberwithin{theorem}{section}

\newcommand{\pif}{purely infinite}

\newcommand{\Id}{\mathrm{id}}
\newcommand{\GL}{\mathrm{GL}}

\newcommand{\spek}{\mathrm{sp}}

\newcommand{\Prim}{\mathrm{Prim}}
\newcommand{\Ideal}{\mathrm{Ideal}}
\newcommand{\Aff}{\mathrm{Aff}}

\newcommand{\RR}{{\rm RR}}

\newcommand{\ep}{{\varepsilon}}

\newcommand{\cO}{{\cal O}}

\newcommand{\cZ}{{\cal Z}}
\newcommand{\C}{{\mathbb C}}

\newcommand{\N}{{\mathbb N}}

\newcommand{\bO}{{\mathbb O}}
\newcommand{\cK}{{\cal K}}

\newcommand{\R}{{\mathbb R}}
\newcommand{\Cs}{{$C^*$-al\-ge\-bra}}

\newcommand{\gange}{\! \cdot \!}

\newcommand{\sh}{{$^*$-ho\-mo\-mor\-phism}}

\newenvironment{proof}[1][Proof:]
{\begin{trivlist}\item[]\textbf{#1} }
{\hbox{}\nobreak\hfill\quad\hbox{$\square$}\end{trivlist}}

\begin{document}

\title{Purely infinite $C^*$-algebras of real rank zero}

\author{Cornel Pasnicu and Mikael R\o rdam}
\date{}
\maketitle
\begin{abstract} \noindent We show that a separable \pif{} \Cs{} is of
  real rank zero if and only if its primitive ideal space has a basis
  consisting of compact-open sets and the natural map
  $K_0(I) \to K_0(I/J)$ is surjective for all closed two-sided ideals
  $J \subset I$ in the \Cs. It follows in particular that if $A$ is any
  separable \Cs{}, then $A \otimes \cO_2$ is of real rank zero if and
  only if the primitive ideal space of 
  $A$ has a basis of compact-open sets, which again
  happens if and only if $A \otimes \cO_2$ has the \emph{ideal property}, also
  known as property (IP). 
\end{abstract}

\section{Introduction} \label{sec:intro}

\noindent The extend to which a \Cs{} contains projections is decisive for
its structure and properties. Abundance of projections can be expressed 
in many ways, several of which were proven to be equivalent by Brown and
Pedersen in \cite{BroPed:realrank}. They refer to \Cs s satisfying
these equivalent conditions as having \emph{real rank zero}, written
$\RR( - ) = 0$, (where the real rank
is a non-commutative notion of dimension). One of these equivalent
conditions states that every non-zero hereditary
sub-algebra has an approximate unit consisting of projections. Real
rank zero is a non-commutative analog of being totally disconnected
(because an abelian \Cs{} $C_0(X)$, where $X$ is a locally compact
Hausdorff space, is of real rank zero if and only if $X$ is
totally disconnected). Another, weaker, condition, that we shall consider
here is the ideal property (denoted (IP)) that projections in the \Cs{}
separate ideals.

The interest in \Cs s of real rank zero comes in parts from the
fact that many \Cs s of interest happen---sometimes surprisingly---to
be of real rank zero, and it comes in parts
from Elliott's classification conjecture which predicts that
separable nuclear \Cs s be classified by some invariant that
includes $K$-theory (and in some special cases nothing more than
$K$-theory!). The Elliott conjecture has a
particularly nice formulation for \Cs s of real rank zero, it has been
verified for a wide class of \Cs s of real rank zero, and the Elliott
conjecture may still hold (in its original form) within this class of \Cs s
(there are counterexamples to Elliott's conjecture in the non-real rank
zero case).

If the Elliott conjecture holds for a certain class of \Cs s, then one
can decide whether a specific \Cs{} in this class is of real rank zero
or not
by looking at its Elliott invariant. In the unital stably finite case,
the Elliott conjecture predicts that a ``nice'' \Cs{} $A$ is of
real rank zero if and only if the image of $K_0(A)$ in $\Aff(T(A))$ is
dense, where $T(A)$ is the simplex of normalized traces on
$A$. This has been verified in \cite{Ror:Z-absorbing} in the case where $A$
in addition is exact and tensorially absorbs the Jiang-Su algebra $\cZ$.
In the presence of some weak divisibility properties on $K_0(A)$, the
condition that $K_0(A)$ has dense image in $\Aff(T(A))$
can be replaced with the weaker condition that projections
in $A$ separate traces on $A$.

In the simple, \pif{} case, where there are no traces, real rank zero is
automatic as shown by Zhang in \cite{Zhang:infsimp}. This
result is here generalized, assuming separability, to the
non-simple case. We are forced to consider obstructions to real rank
zero that do not materialize themselves in the simple case,
including topological properties of the primitive ideal space and
$K$-theoretical obstructions (as explained in the abstract).

The notion of being purely infinite was introduced by Cuntz,
\cite{Cuntz:KOn}, in the
simple case and extended to non-simple \Cs s by Kirchberg and the
second named author in \cite{KirRor:pi} (see Remark~\ref{rm:pi} for
the definition). The study of \pif{} \Cs s was motivated by
Kirchberg's classification of separable, nuclear, (strongly) \pif{}
\Cs s up to stable isomorphism by an ideal related $KK$-theory. This
classification result, although technically and theoretically powerful, is
hard to
apply in practice; however, it has the following beautiful corollary:
Two separable nuclear \Cs s $A$ and $B$ are isomorphic after being tensored
by $\cO_2 \otimes \cK$ if and only if their primitive
ideal spaces are homeomorphic.

Suppose that $A$ is a separable nuclear \Cs{} whose primitive ideal
space has a basis for its topology consisting of compact-open sets.
Then, thanks to a result of Bratteli and Elliott, \cite{BraEll:AF},
there is an AF-algebra $B$ whose primitive ideal space is
homeomorphic to that of $A$. It follows that $A \otimes \cO_2
\otimes \cK \cong B \otimes \cO_2 \otimes \cK$; the latter \Cs{} is
of real rank zero, whence so is the former, whence so is $A \otimes
\cO_2$. In other words, if $A$ is separable and nuclear, then $\RR(A
\otimes \cO_2) = 0$ if and only if the primitive ideal space of $A$
has a basis of compact-open sets. Seeking to give a direct proof of
this result and to drop the nuclearity hypothesis on $A$, we started
the investigations leading to this article.

The paper is divided into three sections. In Section~\ref{sec:ip} we
remind the reader of some of the relevant definitions and concepts, and it
is shown that a \pif{} \Cs{} has property (IP) if and only if its primitive
ideal space has a basis of compact-open sets. Section~\ref{sec:lift}
contains a discussion of the $K$-theoretical obstruction, that we call
\emph{$K_0$-liftable}, to having real rank zero and 
some technical ingredients that are needed for the proof of
our main result, mostly related to lifting properties of projections. The final
Section~\ref{sec:main} contains our main result (formulated in the
abstract) and some corollaries thereof.

Throughout this paper, the symbol $\otimes$ will mean the minimal
tensor product of \Cs{}s; and by an ideal of an arbitrary \Cs{} we
will, unless otherwise specified, mean a closed and two-sided ideal.

\newpage
\section{Purely infinite \Cs s with property (IP)} \label{sec:ip}

\noindent In this section we show, among other things, that a \pif{}
separable \Cs{} has the ideal property if and only if its primitive ideal
space has a basis consisting of compact-open sets. We begin by explaining
the concepts that go into this statement.

\begin{remark}[The ideal property (IP)] \label{rem:ip} A \Cs{} $A$ has
  \emph{the ideal property,} abbreviated (IP), if projections in $A$
  separate ideals in $A$, i.e., whenever $I,J$ are ideals
  in $A$ such that $I \nsubseteq J$, then there is a projection in $I
  \setminus J$.

The ideal property first appeared in Ken Stevens' Ph.D.\ thesis,
where a certain class of (non-simple) \Cs s with the ideal property
were classified by a $K$-theoretical invariant; later the first
named author has studied this concept extensively, see e.g.,
\cite{Pas:ideal_prop_shape} and \cite{Pas:ideal_prop_AH}.
\end{remark}

\begin{remark}[The primitive ideal space] \label{rem:prim}
The primitive ideal space, denoted $\Prim(A)$, of a \Cs{} $A$
is the set of all primitive ideals in $A$ (e.g., kernels of irreducible
representations) equipped with the Jacobsen topology. The
Jacobsen topology is given as follows: if $\mathcal{M} \subseteq \Prim(A)$
and $J \in \Prim(A)$, then
$$J \in \overline{\mathcal{M}} \iff
\bigcap_{I \in {\mathcal{M}}} I \subseteq J.$$
There is a natural lattice isomorphism between the ideal lattice,
denoted $\Ideal(A)$, of $A$ and the lattice, $\bO(\Prim(A))$, of open
subsets of $\Prim(A)$ given as
\begin{eqnarray*}
J \in \Ideal(A) & \leadsto & \{I \in \Prim(A) \colon J \subseteq
I\}^{\mathrm{c}} \in \bO(\Prim(A)), \\
U \in \bO(\Prim(A)) & \leadsto & J = \bigcap_{I \in U^{\mathrm{c}}} I
\in \Ideal(A),
\end{eqnarray*}
(where $U^{\mathrm{c}}$ denotes the complement of $U$). A subset of
$\Prim(A)$ is said to be compact\footnote{Some
authors would rather call such a space \emph{quasi-compact} and reserve the
term ``compact'' for spaces that also are Hausdorff.}
if it has the Heine-Borel property. In the non-Hausdorff setting, compact
sets need not be closed; compactness is preserved under 
forming finite unions, but not under (finite or infinite)
intersections.  

Subsets of $\Prim(A)$ which are both compact and open are, naturally,
called compact-open. An ideal $J$ in $A$ corresponds to a compact-open
subset in $\Prim(A)$ if and only if it has the following property (which is
a direct translation of the Heine-Borel property): Whenever $\{J_\alpha\}$
is an increasing net of ideals in $A$ such that $J =
\overline{\bigcup_\alpha J_\alpha}$, then $J = J_\alpha$ for some
$\alpha$. We shall often---sloppily---refer to such ideals as
\emph{compact ideals}. 

We are particularly interested in the case where $\Prim(A)$ has a
basis (for its topology) consisting of compact-open sets. When
$\Prim(A)$ is locally compact and Hausdorff this is the case precisely
when $\Prim(A)$ is
totally disconnected (all connected components are singletons). In
general, $\Prim(A)$ has a basis of compact-open sets if and only if
every (non-empty) open subset is the union of an increasing
net of compact-open sets, or, equivalently, if and
only if every ideal $J$ in $A$ is equal to $\overline{\bigcup_\alpha
  J_\alpha}$ for some increasing net $\{J_\alpha\}_\alpha$ of compact ideals.

If $\Prim(A)$ is finite, which happens precisely when $\Ideal(A)$ is
finite, then all subsets are compact, whence $\Prim(A)$ has a basis of
compact-open sets. The space $\Prim(A)$ is totally disconnected in
this case if and only if it is Hausdorff; or, equivalently, if and
only if $A$ is the direct sum of finitely many simple \Cs s. See also
Example~\ref{ex:finite}. 

If $A$ is a separable \Cs{}, then $\Prim(A)$ is a locally compact
second countable T$_0$-space in which every (closed) prime\footnote{A set $F$
  is called prime if whenever $G$ and $H$ are closed sets with $F = G
  \cup H$, then $F=G$ or $F=H$.} subset is the closure of a
point. Conversely, if $X$ is a space with these properties, and if $X$
has a basis for its topology consisting of compact-open sets, then $X$
is homeomorphic to $\Prim(A)$ for some separable AF-algebra $A$, as shown by
Bratteli and Elliott in \cite{BraEll:AF}.
\end{remark}

\noindent We shall need the following (probably well-known) easy lemma:

\begin{lemma} \label{lm:totdisc}
Let $A$ be a \Cs, let $I, I_1, I_2$ be ideals in $A$, and
let $\pi \colon A \to A/I$ be the quotient mapping.
\begin{enumerate}
\item If $I_1$ and $I_2$ are compact, then so is $I_1+I_2$.
\item If $I$ is compact and if $J$ is a compact ideal in $A/I$, then
  $\pi^{-1}(J)$ is compact.
\end{enumerate}
\end{lemma} 

\begin{proof} (i). The union of two compact sets is again compact (also in
  a T$_0$-space).

(ii). Let $\{K_\alpha\}_\alpha$ be an arbitrary upwards directed family of
ideals in $A$ such that $\bigcup_\alpha K_\alpha$ is dense
in $\pi^{-1}(J)$. Then $J = \overline{\bigcup_\alpha \pi(K_\alpha)}$, whence
  $J = \pi(K_{\alpha_1})$ for some $\alpha_1$. As $I$ is contained in
  $\pi^{-1}(J)$, it is equal to the closure of $\bigcup_\alpha (I \cap
  K_\alpha)$, whence $I = I \cap K_{\alpha_2}$ for some $\alpha_2$. It now
  follows that $\pi^{-1}(J) = K_\alpha$ whenever $\alpha$ is greater than
  or equal to both $\alpha_1$ and $\alpha_2$. 
\end{proof}

\noindent We shall show later (in Corollary~\ref{cor:extension}) 
that the class of \Cs s,
for which the primitive ideal space has a basis of compact-open sets, is
closed under extensions.

\begin{remark}[Scaling elements] \label{rem:scaling}
Scaling elements were introduced by Blackadar and Cuntz in
\cite{BlaCuntz:infproj} as a mean to show the existence of projections in
simple \Cs s that admit no dimension function. An element $x$ in a \Cs{}
$A$ is called a scaling element if $x$ is a contraction and
$x^*x$ is a unit for $xx^*$, ie., if $x^*xxx^* = xx^*$. Blackadar and Cuntz
remark that if $x$ is a scaling
element, then $v = x + (1-x^*x)^{1/2}$ is an isometry in the unitization of
$A$, whence $p = 1-vv^*$ is a projection in $A$. Moreover, if $a$ is a
positive element in $A$ such that $x^*xa = a$ and  $xx^*a=0$, then
$pa=a$. In this way we get a ``lower bound'' on the projection $p$.
\end{remark}

\begin{remark}[Cuntz' comparison theory] \label{rm:comparison}
We recall briefly the notion of comparison of positive elements in a
\Cs{} $A$, due to Cuntz, \cite{Cuntz:dimension}.
Given $a,b \in A^+$, write $a \precsim b$ if
for all $\ep >0$ there is $x \in A$ such that $\|x^*bx-a\| < \ep$. Let
$(a-\ep)_+$ denote the element obtained by applying the function $t
\mapsto \max\{t-\ep,0\}$ to $a$. It is shown in \cite{Ror:UHFII} that
if $a,b$ are positive elements in $A$ and if $\ep >0$, then 
$\|a-b\| < \ep$ implies $(a-\ep)_+ \precsim b$; and $a \precsim b$ and
$a_0 \in \overline{(a-\ep)_+A(a-\ep)_+}$ implies that $a_0 =
x^*bx$ for some $x \in A$.

We shall also need the following fact: If $a \precsim b$ and $\ep > 0$,
then there exists a contraction $z \in A$ such that
$z^*z(a-\ep)_+ = (a-\ep)_+$ and $zz^* \in \overline{bAb}$. Indeed, there is
a positive contraction $e$ in 
$\overline{(a-\ep/2)_+A(a-\ep/2)_+}$ 
such that $e(a-\ep)_+ = (a-\ep)_+$, and by the result mentioned above 
there is $x \in A$ such that
$e = x^*bx$. The element $z = b^{1/2}x$ is now as desired.
\end{remark}

\begin{remark}[Purely infinite \Cs s] \label{rm:pi}
A (possibly non-simple) \Cs{} $A$ is said to be \emph{purely infinite}
if $A$ has no character (or, equivalently, no abelian quotients) and
if
$$\forall a,b \in A^+ : a \in \overline{AbA} \iff a \precsim b,$$
where $\overline{AbA}$ denotes the ideal in $A$
generated by the element $b$. Observe that the implication
``$\Leftarrow$'' above is trivial and holds for all \Cs s.

It is shown in \cite{KirRor:pi} that any positive element $a$ in a
\pif{} \Cs{} is properly infinite (meaning that $a \oplus a \precsim
a \oplus 0$ in $M_2(A)$); and in particular, all (non-zero)
projections in a \pif{} \Cs{} are properly infinite (in the standard
sense: $p \in A$ is properly infinite if there are projections
$p_1,p_2 \in A$ such that $p_j \le p$, $p_1 \perp p_2$, and $p_1 \sim
p_2 \sim p$).

It is also proved in \cite{KirRor:pi} that $A \otimes \cO_\infty$ and
$A \otimes \cO_2$ are \pif{} for all \Cs s $A$, and hence that $A
\otimes B$ is \pif{} whenever $B$ is a Kirchberg algebra\footnote{A
  simple, separable, nuclear, \pif{} \Cs.} (because these satisfy $B
\cong B \otimes \cO_\infty$, see \cite{KirPhi:classI}).
\end{remark}

\begin{proposition} \label{prop:compact}
Let $I$ be an ideal in a separable \pif{} \Cs{}
  $A$. Then the following conditions are equivalent:
\begin{enumerate}
\item $I$ corresponds to a compact-open subset of $\Prim(A)$, i.e., $I$ is
  compact. 
\item $I$ is generated by a single projection in $A$.
\item $I$ is generated by a finite family of projections in $A$.
\end{enumerate}
\end{proposition}

\begin{proof}
(iii) $\Rightarrow$ (i). Suppose that $I$ is
generated (as an ideal) by the projections $p_1, \dots, p_n$,
and suppose that $\{I_\alpha\}_\alpha$
is an increasing net of ideals in $A$ such that
$\bigcup_\alpha I_\alpha$ is a dense (algebraic) ideal in $I$. Then
$\bigcup_\alpha 
I_\alpha$ contains the projections $p_1, \dots, p_n$ (because it contains
the Pedersen ideal\footnote{This is the
smallest dense algebraic two-sided ideal in the \Cs.} 
of $I$, and the Pedersen ideal of a \Cs{} contains
all projections of the  \Cs{}). It follows that $p_1, \dots, p_n$ belong
to $I_\alpha$ for some $\alpha$, whence $I = I_\alpha$.

(i) $\Rightarrow$ (ii). By separability of $A$ (and hence
of $I$), $I$ contains a strictly positive element, and is hence generated
(as an ideal) by a single positive element $a$. For each
$\ep \ge 0$ let $I_\ep$ be the ideal in $A$ generated by
$(a-\ep)_+$. Then $I = \overline{\bigcup_{\ep >0} I_\ep}$, so by assumption
(and Remark~\ref{rem:prim}), $I = I_{\ep_0}$ for some $\ep_0 > 0$. It
follows in particular that $a \precsim (a-\ep_0)_+$, cf.\
Remark~\ref{rm:pi}.

Choose $\ep_1$ such that $0 < \ep_1 < \ep_0$. As $A$ is \pif, all
its positive elements, and in particular $(a-\ep_1)_+$, are properly
infinite (see \cite[Definition~3.2]{KirRor:pi}). Use
\cite[Proposition~3.3]{KirRor:pi} (and
\cite[Lemma~2.5~(i)]{KirRor:pi}) to find mutually orthogonal
positive elements $b_1, b_2$ in $\overline{(a-\ep_1)_+I(a-\ep_1)_+}$
such that $(a-\ep_0)_+ \precsim b_1$ and $(a-\ep_0)_+ \precsim b_2$.
Then $a \precsim b_1$ (a relation that also holds relatively to $I$) and $a
\precsim b_2$ (whence $b_2$ is full in $I$). 

By Remark~\ref{rm:comparison} there is $x \in I$ such that
$x^*x(a-\ep_1)_+ = (a-\ep_1)_+$ and $xx^*$ belongs to $\overline{b_1
I b_1} \subseteq \overline{(a-\ep_1)_+A(a-\ep_1)_+}$. We conclude
that $x$ is a scaling element which satisfies $x^*x b_2 = b_2$ and
$xx^* b_2 =0$. By the result of Blackadar and Cuntz mentioned in
Remark~\ref{rem:scaling} above there is a projection $p \in I$ such
that $pb_2 = b_2$. As $b_2$ is full in $I$, so is $p$.

(ii) $\Rightarrow$ (iii) is trivial.
\end{proof}

\noindent The implications (ii) $\Rightarrow$ (iii) $\Rightarrow$
(i) do not require separability of $A$.

The corollary below follows immediately from
Proposition~\ref{prop:compact} (and from Remark~\ref{rem:prim}).

\begin{corollary} \label{cor:ip}
Let $A$ be a separable \pif{} \Cs{} where $\Prim(A)$ has a basis of
compact-open sets. Then any ideal in $A$ is either generated by
a single projection or is the closure of the union of an increasing
net of ideals each of which is generated by a single projection.

Converserly, if $A$ is any \Cs{} (not necessarily separable or \pif), and
if any ideal in $A$ either is generated by a single projection or is
the closure of the union of an increasing net of ideals with this property,
then $\Prim(A)$ has a basis for its topology consisting of compact-open
sets. 
\end{corollary}

\begin{lemma} \label{lm:her} Let $A$ be a \pif{} \Cs{} and let $B$ be a
  hereditary sub-\Cs{} of $A$. Then each projection in
  $\overline{ABA}$, the ideal in $A$ generated by
  $B$, is equivalent to a projection in $B$.
\end{lemma}

\begin{proof} Let $p$ be a projection in $\overline{ABA}$. The family
  of ideals in $A$
  generated by a single positive element in $B$ is upwards directed
  (if $I_1$ is generated by $b_1 \in B^+$ and $I_2$ is generated by $b_2
  \in B^+$, then $I_1+I_2$ is generated by $b_1 +b_2 \in B^+$). The union
  of these ideals is dense in $I$ and therefore contains $p$. It follows
  that $p$ belongs to $\overline{AbA}$ for some $b \in B^+$. As $A$ is
  \pif, $p \precsim b$, whence $p = z^*bz$ for some $z
  \in A$ (because $p$ is a projection, cf.\
  \cite[Proposition~2.7]{Ror:UHFII}). Put $v = b^{1/2}z$. Then
  $v^*v = p$, and $vv^* \in B$ is therefore a projection which is
  equivalent to $p$.
\end{proof}

\begin{proposition} \label{prop:ipher}
Any hereditary sub-\Cs{} of a \pif{} \Cs{} with property (IP) again
has property (IP).
\end{proposition}

\begin{proof} Let $A_0$ be a hereditary sub-\Cs{} of a  \pif{} \Cs{} $A$
with property (IP), and let $I_0$ and $J_0$ be
ideals in $A_0$ with $I_0 \nsubseteq J_0$. Let $I$ and $J$ be the
ideals in $A$
generated by $I_0$ and $J_0$, respectively. Then $I \nsubseteq J$, and
so, by assumption, there is a projection $p \in I \setminus J$. By
Lemma~\ref{lm:her}, $p$ is equivalent to a projection $p' \in I_0$;
and $p'$ does not belong $J$, and hence not to $J_0$.
\end{proof}

\noindent Condition (iii) below was considered by Brown and Pedersen in
\cite[Theorem~3.9 and Discussion~3.10]{BroPed:idealstructure} and was there
given the name \emph{purely properly infinite}. Brown and Pedersen noted
that purely properly infinite \Cs s are purely infinite (in the sense
discussed in Remark~\ref{rm:pi}). Brown kindly informed us that this
property is equivalent with properties (i) and (ii) below. We thank Larry
Brown for allowing us to include this statement here. 

\begin{proposition} \label{prop:IP}
The following four conditions are equivalent for any separable \Cs{} $A$.
\begin{enumerate}
\item $A$ is \pif{} and $\Prim(A)$ has a basis for its topology consisting of
  compact-open sets.
\item $A$ is \pif{} and has property (IP).
\item Any non-zero hereditary sub-\Cs s of $A$ is generated as an ideal by
  its properly infinite projections. 
\item Every non-zero hereditary sub-\Cs{} in any quotient of $A$ contains
  an infinite projection. 
\end{enumerate}
The implications {\rm{(i)}} $\Leftarrow$ {\rm{(ii)}} $\Leftrightarrow$
{\rm{(iii)}} 
$\Leftrightarrow$ {\rm{(iv)}} hold also when $A$ is non-separable. 
\end{proposition}

\begin{proof} Separability is assumed only in the proof of ``(i)
  $\Rightarrow$ (ii)''. 

(ii) $\Rightarrow$ (i). Let $I$ be an ideal in $A$. Then $I$ is generated
by its projections (because $A$ has property (IP)). Let $\Lambda$ be the
net of finite subsets 
  of the set of projections in $I$, and, for each $\alpha \in \Lambda$, let
  $I_\alpha$ be the ideal in $A$ generated by the
  projections in the finite set $\alpha$. Then $I_\alpha$ is compact (by
  Proposition~\ref{prop:compact}), and $\bigcup_{\alpha \in \Lambda}
  I_\alpha$ is dense in $I$. This shows that $\Prim(A)$ has a basis of
  compact-open sets, cf.\ Remark~\ref{rem:prim}. 

(i) $\Rightarrow$ (ii). Suppose that (i) holds, and let
  $I,J$ be ideals in $A$ such that $I \nsubseteq J$.
From Corollary~\ref{cor:ip} there is an increasing net of ideals $I_\alpha$ in
$A$ each generated by a 
single projection, say $p_\alpha$, such that
  $\bigcup_\alpha I_\alpha$ is dense in $I$. Now, $I_\alpha \nsubseteq J$
  for some $\alpha$, and so the projection $p_\alpha$ belongs to
$I \setminus J$.

(ii) $\Rightarrow$ (iii). Every non-zero projection in a \pif{} \Cs{} is
properly infinite (see Remark~\ref{rm:pi} or \cite[Theorem~4.16]{KirRor:pi})
and so it suffices to show that any hereditary sub-\Cs{} of $A$ has
property (IP); but this follows from Proposition~\ref{prop:ipher} and the
assumption that $A$ is \pif{} and has property (IP).   

(iii) $\Rightarrow$ (iv). Let $I$ be an ideal in $A$, and
let $B$ be a non-zero hereditary sub-\Cs{} of $A/I$. Let $\pi \colon A \to
A/I$ denote the quotient mapping. By (iii) and Lemma~\ref{lm:her} there is
a properly infinite 
projection $p$ in $\pi^{-1}(B) \setminus I$; and so $\pi(p)$ is a non-zero
properly infinite (and hence infinite) projection in $B$.

(iv) $\Rightarrow$ (i). It follows from \cite[Proposition~4.7]{KirRor:pi}
that $A$ is \pif{}. We must show that $\Prim(A)$ has a basis of
compact-open sets. We use the equivalent formulation given in
Remark~\ref{rem:prim}, see also Corollary~\ref{cor:ip}. 

Let $I$ be an ideal in $A$, and let $\{I_\alpha\}$ be the family of
all compact ideals contained in $I$. Then
$\{I_\alpha\}_\alpha$ is upwards directed (by
Lemma~\ref{lm:totdisc}~(i)). Let $I_0$ be the closure of the union of the
ideals $I_\alpha$. We must show that $I_0 = I$. Suppose, to reach a
contradiction, that $I_0 \subset
I$. Then, by (iv), $I/I_0$ contains a non-zero projection $p$. The
projection $p$ lifts to a projection $q$ in $I/I_\alpha$ for some
$\alpha$ (by semiprojectivity of the \Cs{} $\C$, see also the proof of
Lemma~\ref{lm:2} below). Let $J$ be the ideal in $I/I_\alpha$
generated by the projection $q$. Then $J$ 
is compact, whence so is its pre-image $I' \subseteq I$ under the quotient
mapping $I \to I/I_\alpha$, cf.\ Lemma~\ref{lm:totdisc}~(ii). As the image
of $I'$ under the quotient mapping $I \to I/I_0$ contains the projection
$p$ we conclude that $I'$ is not contained in $I_0$, which is in
contradiction with the construction of $I_0$.   
\end{proof}

\noindent Property (i) in the lemma below is pretty close to saying that
the hereditary sub-\Cs{} $\overline{aAa}$ has an approximate unit
consisting of projections, and hence that $A$ is of real rank zero. In
fact, if $A$ has stable rank one (which by the way never can happen when
$A$ is purely infinite and not stably projectionless!), then property (i)
below would have implied that $A$ has real rank zero. In the absence of
stable rank one we get real rank zero from condition (i) below if a
$K$-theoretical condition, discussed in the next section, is satisfied.

\begin{lemma} \label{lm:proj2}
Let $A$ be a \pif{} \Cs{} with property (IP).
\begin{enumerate}
\item For each positive
element $a \in A$ and for each $\ep > 0$, there is a projection $p \in
\overline{aAa}$ such that $(a-\ep)_+ \precsim p$.
\item For each element $x \in A$ and for each $\ep > 0$, there is a
  projection $p \in A$ and an element $y \in A$ such that $\|x-y\| \le
  \ep$ and $y \in \overline{ApA}$.
\end{enumerate}
\end{lemma}

\begin{proof} (i). The hereditary \Cs{} $\overline{aAa}$ is \pif{} and
has property (IP) (by Lem\-ma~\ref{lm:her}). We can therefore apply
Corollary~\ref{cor:ip} to  $\overline{aAa}$ to obtain an increasing
net $\{I_\alpha\}_\alpha$ of ideals in
$\overline{aAa}$ each generated by a single projection such that
$\bigcup_\alpha I_\alpha$ is a dense algebraic ideal in
$\overline{aAa}$. It follows 
that $(a-\ep)_+$ belongs to $\bigcup_\alpha I_\alpha$, and hence to
$I_\alpha$ for some $\alpha$. Let $p$ be a projection that generates
the ideal $I_\alpha$. Then $(a-\ep)_+ \precsim p$, because $(a-\ep)_+$
belongs to the ideal generated by $p$.

(ii). Write $x = v|x|$ with $v$ a partial isometry in $A^{**}$, and
put $y = v(|x|-\ep)_+ \in A$. Then $\|x-y\| \le \ep$ and $|y| =
(|x|-\ep)_+$. Use (i) to find a projection $p$ in $A$ such that $|y|
\precsim p$. Then $|y|$, and hence also $y$, belong to $\overline{ApA}$.
\end{proof}

\noindent We continue this section with a general result on \Cs s
(not necessarily purely infinite) with property (IP) that is
relevant for the discussion in Section~\ref{sec:lift}.

\begin{proposition} \label{prop:stable}
Any separable stable \Cs{} with property (IP) has an approximate unit
consisting of projections.
\end{proposition}

\begin{proof} If $A$ is a separable stable \Cs{} containing a full
  projection $p$, then $A$ is isomorphic to $pAp \otimes \cK$ by Brown's
  theorem; and so in particular $A$ has an approximate unit consisting of
  projections.

Suppose that $A$ is separable, stable and with property
  (IP). Then $A = \overline{\bigcup_\alpha A_\alpha}$ for some increasing
  net $\{A_\alpha\}_\alpha$ of ideals in $A$ each of which
  is generated by a finite set of projections, cf.\ the proof of ``(i)
  $\Rightarrow$ (ii)'' in Proposition~\ref{prop:IP}. We claim that each
  $A_\alpha$ is in fact generated by a single projection. Indeed, suppose
  that $A_\alpha$ is generated as an ideal by the projections $p_1, p_2,
  \dots, p_n$; then $p_1 \oplus p_2 \oplus \cdots \oplus p_n$ is equivalent
  to (or equal to) a projection $p \in A_\alpha$, because $A_\alpha$ is
  stable (being an ideal in a stable \Cs{}). It follows that $A_\alpha$ is
  generated by the projection $p$. By the first part of the proof,
  $A_\alpha$ has an approximate
  unit consisting of projections. As this holds for all $\alpha$ we
  conclude  that also $A$ has an approximate unit consisting of projections.
\end{proof}

\noindent Proposition~\ref{prop:3pi} below was shown in \cite{KirRor:pi2} 
by Kirchberg and the second named author for \Cs s of the real
rank zero. We extend here this result to the broader class of \Cs s
with property (IP). We refer to  \cite{KirRor:pi2} for the definitions
of being strongly, respectively, weakly purely infinite. 

\begin{proposition} \label{prop:3pi}
    Let $A$ be a \Cs{} with property (IP).  The following are
    equivalent:
    \begin{enumerate}
        \item $A$ is purely infinite.
        \item $A$ is strongly purely infinite.
        \item $A$ is weakly purely infinite.
    \end{enumerate}
\end{proposition}

\begin{proof}
(ii) $\Rightarrow$ (i) $\Rightarrow$ (iii) are (trivially) true for all \Cs
s $A$ (see \cite[Theorem~9.1]{KirRor:pi2}). 

    (i) $\Rightarrow$ (ii).  It follows from Lemma~\ref{lm:proj2} and from
    \cite[Remark~6.2]{KirRor:pi2} (see also the proof of
    \cite[Proposition~6.3]{KirRor:pi2}) that any \Cs{} with property (IP) has
    the \emph{locally central
    decomposition property}; and \cite[Theorem~6.8]{KirRor:pi2} says that any
  \pif{} \Cs{} with the locally central decomposition property is strongly
  \pif. 

    (iii) $\Rightarrow$ (i).  Assume that $A$ is weakly purely
    infinite.  By the comment following
    \cite[Proposition~4.18]{KirRor:pi2}, $A$ is purely infinite if
    every quotient of $A$ has the property (SP) (i.e., each non-zero
    hereditary sub-\Cs{} contains a non-zero projection).  
    Both property (IP) and weak pure infiniteness pass to
    quotients, cf.\ \cite[Proposition~4.5]{KirRor:pi2}, so it will be enough
    to prove that any non-zero hereditary sub-\Cs{}
    $B$ of $A$ contains a non-zero projection.  

As $A$ is weakly \pif{}, it is pi-$n$ for some natural number $n$ (see
\cite[Definition~4.3]{KirRor:pi2}). By the Glimm lemma (see
\cite[Proposition~4.10]{KirRor:pi}) there is a non-zero
   \sh{} from $M_n(C_0((0,1]))$ into $B$.  So we get
    non-zero pairwise equivalent and orthogonal positive elements $e_1,
    \dots, e_n$ in $B$. The ideal in $A$ generated by
    $e_1$ contains a non-zero projection $p$. As $A$ is assumed to be
    pi-$n$ we can use 
    \cite[Lemma~4.7]{KirRor:pi2} to conclude that $p \precsim e_1 \otimes 1_n$;
    and as $e_1 \otimes 1_n \precsim e_1 + e_2 + \cdots + e_n =: b \in B$
    (see \cite[Lemma~2.8]{KirRor:pi}) it follows from 
    \cite[Proposition~2.7~(iii)]{KirRor:pi} that $p$ is equivalent
    to a (necessarily non-zero) projection $q$ in $\overline{bAb} \subseteq
    B$. (It has been used twice above that $p \precsim (1-\ep)p =
    (p-\ep)_+$ when $p$ is a projection and $0 \le \ep < 1$.) 
\end{proof}

\section{Lifting projections}  \label{sec:lift}

\noindent We consider here when projections in a quotient of a \pif{} \Cs{}
lift to the
\Cs{} itself. We begin with a discussion of a $K$-theoretical
obstruction to lifting projections: 

\begin{definition} \label{def:K0onto}
A \Cs{} $A$ is said to be
\emph{$K_0$-liftable} if for every pair of ideals
$I \subset J$ in $A$, the extension
$$\xymatrix{0 \ar[r] & I \ar[r]^-\iota & J \ar[r]^-\pi & J/I \ar[r] & 0}$$
has the property that $K_0(\pi) \colon K_0(J) \to K_0(J/I)$ is
surjective (or, equivalently, that the index map $\delta \colon
K_0(J/I) \to K_1(I)$ is zero, or, equivalently, if the induced map
$K_1(\iota) \colon K_1(I) \to K_1(J)$ is injective).
\end{definition}

\noindent As pointed out to us by Larry Brown, it suffices to check
$K_0$-liftability for $J=A$ (i.e., $A$ is $K_0$-liftable if and only if
the induced map $K_0(A) \to K_0(A/I)$ is onto for every ideal $I$ in $A$),
because if $K_1(I) \to K_1(A)$ is injective, then so is $K_1(I) \to K_1(J)$
whenever $I \subseteq J \subseteq A$.

Every simple \Cs{} is automatically
$K_0$-liftable (there are no non-trivial sequences $0 \to I \to J
\to J/I \to 0$ for ideals $I \subset J$ in a simple \Cs{}).

The property real rank zero passes from a \Cs{} to its ideals (cf.\
Brown and Pedersen, \cite{BroPed:realrank}), and in the
same paper it is shown that the map $K_0(A) \to K_0(A/I)$ is onto
whenever $A$ is a \Cs{} of real rank zero and $I$ is an ideal in
$A$. Hence all \Cs s of real rank zero are $K_0$-liftable.

Being $K_0$-liftable passes to hereditary sub-\Cs s:

\begin{lemma} \label{lm:K0onto}
Any hereditary sub-\Cs{} of a separable $K_0$-liftable \Cs{} is
again $K_0$-liftable.
\end{lemma}

\begin{proof}
    Let $A$ be a separable $K_0$-liftable \Cs, and let
    $A_0$ be a hereditary sub-\Cs{} of $A$. Let $I_0 \subset J_0$ be
    ideals in $A_0$, and let $I \subset J$ be the
    ideals in $A$ generated by $I_0$ and $J_0$, respectively.  

Then $J_0$ is a full hereditary sub-\Cs{} of $J$, and (the image in
$J/I$ of) $J_0/I_0$ is a full hereditary sub-\Cs{} in $J/I$. 
The commutative diagram
$$\xymatrix{J_{0} \ar[d] \ar[r] & J_{0}/I_{0} \ar[d] \\
J \ar[r] & J/I }$$
induces a commutative diagram of $K_0$-groups
$$\xymatrix{K_0(J_0) \ar[d] \ar[r]& K_0(J_0/I_{0})\ar[d]
\\ K_0(J) \ar[r] & K_0(J/I),}$$
where the vertical maps are isomorphisms (by stability of $K_0$ and by
Brown's theorem) and
the lower horizontal map is surjective by assumption.
Hence the upper horizontal map $K_0(J_0) \to K_0(J_0/I_{0})$ is
surjective.
\end{proof}

\noindent The next lemma expresses when an extension of two
$K_0$-liftable \Cs s is $K_0$-liftable:

\begin{lemma} \label{lm:K0lift-extension}
Let 
$$\xymatrix{0 \ar[r] & I \ar[r] & A \ar[r]^-\pi & B \ar[r] & 0}$$
be a short-exact sequence of \Cs s. Then $A$ is $K_0$-liftable if and
only if $I$ and $B$ are $K_0$-liftable and the induced  map $K_0(A)
\to K_0(B)$ is onto. 
\end{lemma}

\begin{proof} 
``If''.  We use the remark below Definition~\ref{def:K0onto} whereby it
suffices to show that the map $K_0(A) \to K_0(A/J)$ is onto whenever $J$ is
an ideal in $A$. To this end, consider the diagram of \Cs s with exact
rows and columns: 
$$\xymatrix@C-0.5pc@R-0.5pc{ & 0 \ar[d] & 0 \ar[d] & 0 \ar[d] & \\
0 \ar[r] & I \cap J \ar[d] \ar[r] & I  \ar[d] \ar[r] & I / I
\cap J \ar[d] \ar[r] & 0 \\
0 \ar[r] & J \ar[d] \ar[r] & A \ar[d] \ar[r] & A / J \ar[d] \ar[r] & 0 \\
0 \ar[r] & \pi(J) \ar[d] \ar[r] & B \ar[d] \ar[r] & B / \pi(J)
\ar[d] \ar[r] & 0 \\ & 0 & 0 & 0 &}
$$
that induces the following diagram at the level of $K_0$:
$$
\xymatrix{K_0(I) \ar[r]^-{*} \ar[d] & K_0(I  / I \cap J) \ar[d]
\\ K_0(A) \ar[d]_-{*} \ar[r] & K_0(A/J) \ar[d] \\
K_0(B) \ar[r]^-{*} & K_0(B/\pi(J))}
$$
where the vertical sequences are exact and the maps marked with an asterisk
are surjective (by our assumptions). A
standard diagram chase shows that the map $K_0(A) \to K_0(A/J)$ is
surjective.

``Only if''. If $A$ is $K_0$-liftable, then clearly so
  is $I$, and $K_0(A) \to K_0(B)$ is onto. We proceed to prove that $B$
  is $K_0$-liftable. Let $J \subset L$ be ideals in $B$, and consider
  the commuting diagram
$$\xymatrix{\pi^{-1}(L) \ar[r] \ar[d] & L \ar[d] \\
\pi^{-1}(L)/\pi^{-1}(J) \ar[r]^-{\cong} & L/J,}$$
that induces the commuting diagram
$$\xymatrix{K_0(\pi^{-1}(L)) \ar[r] \ar[d] & K_0(L) \ar[d] \\
K_0\big(\pi^{-1}(L)/\pi^{-1}(J)\big) \ar[r]^-{\cong} & K_0(L/J).}$$
The left-most vertical map is onto by $K_0$-liftability of $A$,
which implies surjectivity of the right-most vertical map. 
\end{proof}

\noindent We proceed to describe when certain tensor products are
$K_0$-liftable

\begin{lemma} \label{lm:K0lift-tensor}
The tensor product $A \otimes \cO_2$ is $K_0$-liftable for every
\Cs{} $A$; and the tensor product $A \otimes \cO_\infty$ is
$K_0$-liftable if and only if $A$ itself is $K_0$-liftable. 
\end{lemma}

\begin{proof} If $D$ is a simple nuclear \Cs{}, then the
  mapping $I \mapsto I \otimes D$ defines a lattice isomorphism from
  $\Ideal(A)$ onto $\Ideal(A \otimes D)$ (surjectivity follows from a
  theorem of Blackadar, \cite{Bla:tensor}, see also
  \cite[Proposition~2.16]{BlanKir:pi3}). Moreover, by Blackadar's
  theorem or by exactness of $D$, if $I \subset J$
  are ideals in $A$, then $(J \otimes D)/(I \otimes D)$ is isomorphic
  to $(J/I) \otimes D$. Hence, to prove $K_0$-liftability
  of $A \otimes D$ it suffices to show that the induced map $K_0(J
  \otimes D) \to K_0((J/I) \otimes D)$ is surjective, or, equivalently,
  that the index map $K_0((J/I) \otimes D) \to K_1(I \otimes D)$ is
  zero. The latter holds for
  all \Cs s $A$ if $D = \cO_2$ because $K_1(I \otimes \cO_2) = 0$. 

To prove the last statement, consider the commutative diagram
$$\xymatrix{J \ar[r] \ar[d] & J/I \ar[d] \\ J \otimes \cO_\infty \ar[r]
  & (J/I) \otimes \cO_\infty,}
$$
where the vertical maps are defined by $x \mapsto x \otimes 1$. It
follows from the K\"unneth theorem that the vertical maps
above induce isomorphisms at the level of $K_0$. It is now clear that
$A \otimes \cO_\infty$ is $K_0$-liftable if and only if $A$ is
$K_0$-liftable. 
\end{proof}

\noindent We now proceed with the projection lifting results. We need
a sequence of lemmas.

\begin{lemma} \label{lm:inv}
Let $A$ be a \Cs, let $x$ be an element in $A$, and suppose that there
is a positive element $e$ in $A$ such that $x^*x \precsim e$, and $x^*x$ and
$xx^*$ are orthogonal to $e$.
Then $x$ belongs to the closure of the invertible elements,
$\GL(\widetilde{A})$, in the unitization $\widetilde{A}$ of $A$.
\end{lemma}

\begin{proof} Let $\ep>0$ be given. By the assumption that
$|x|^2= x^*x \precsim e$ and by Remark~\ref{rm:pi} we obtain a
contraction $z \in A$ such that
$$(|x|-\ep)_+z^*z = (|x|-\ep)_+, \qquad zz^* \in \overline{eAe},
\qquad zz^* \perp z^*z.$$
Now, $a = z+z^*$ is a self-adjoint contraction in $A$ and
$u = a + i\sqrt{1-a^2}$ is a unitary element in $\widetilde{A}$.

Write $x = v|x|$ with $v$ a partial isometry in $A^{**}$, and put
$x_\ep = v(|x|-\ep)_+ \in A$. Then $\|x-x_\ep\| \le \ep$,
$$x_\ep u = v(|x|-\ep)_+ u =  v(|x|-\ep)_+ a =  v(|x|-\ep)_+ z^*,$$
$z^*x_\ep=0$, and so
$$(x_\ep u)^2 = v(|x|-\ep)_+ z^* x_\ep u =0.$$
It follows that $x_\ep  + \lambda u^* = (x_\ep u + \lambda 1)u^*$ is
invertible in $\widetilde{A}$ for all $\lambda \ne 0$, whence $x_\ep$
belongs to the closure of $\GL(\widetilde{A})$. As $\ep > 0$ was
arbitrary, the lemma is proved.
\end{proof}

\begin{lemma} \label{lm:apprunit}
Let $A$ be a \Cs{}. Let $x$ be an element in $A$ and let
$e$ be a properly infinite projection in $A$ such that $x^*x$ is orthogonal
  to $e$ and $x^*x \precsim e$. Then, for each $\ep >0$, there is a
  projection $p \in A$ such that $\|x-xp\| \le \ep$.
\end{lemma}

\begin{proof} Because $e$ is properly infinite (cf.\
  Remark~\ref{rm:pi}) there is a subprojection $e_0$ of $e$ such that
  $e \precsim e_0$ and $e \precsim e-e_0$.
As $|x|^2 = x^*x \precsim e \precsim e_0$ there is $z \in A$ with
$z=e_0z$ and
$$(|x|-\ep/2)_+ = z^*e_0z =z^*z$$
(see Remark~\ref{rm:comparison}).
As $zz^*$ and $z^*z$ both are orthogonal to the projection
$e-e_0$, and $z^*z \precsim e_0 \precsim e-e_0$,  we 
conclude from Lemma~\ref{lm:inv} that $z$ belongs to the closure of
$\GL(\widetilde{A})$. By \cite{Ror:unitary} there is a unitary $u$ in
$\widetilde{A}$ such that
$$ u(|x|-\ep)_+u^* =  u(z^*z-\ep/2)_+u^* = (zz^*-\ep/2)_+ \in e_0Ae_0.$$
The projection $p = u^*e_0u \in A$ thus satisfies $(|x|-\ep)_+ p=
(|x|-\ep)_+$, which entails that
$$\|x(1-p)\| = \| |x|(1-p)\| \le \ep.$$
\end{proof}

\noindent The lemma below and its proof are similar to
\cite[Lemma~3.13]{BroPed:realrank} and its proof.

\begin{lemma} \label{lm:BP}
Let $A$ be a \pif{} \Cs, let $I$ be an ideal
  in $A$, and let $B$ be a hereditary sub-\Cs{} of $A$. Assume that
  $I$ has property (IP). Let $p$ be a
  projection in $B+I$ and assume that $B \cap pAp$ is full in $A$. Then
  there is a projection $q \in B$ such that $p-q$ belongs to $I$.
\end{lemma}

\begin{proof}
Write $p = b + x$, with $b$ a self-adjoint element in $B$ and $x$ a
self-adjoint element in $I$. Take $\ep >0$ such that $2\|b\|\ep + \ep^2 <
1/2$. By Lemma~\ref{lm:proj2} we can find an element
$y \in I$ and a projection $f \in I$ such that $\|x-y\| < \ep/2$
and such that $y$ belongs to the ideal $I_0$ in $I$ generated by $f$.
By assumption, $B \cap pAp$ is full in $A$, so $f$ is equivalent to a
projection $g \in B \cap pAp$ (by Lemma~\ref{lm:her}). Put
$$b_1 =  (1-g)b(1-g) + g \in B, \quad x_1 = (1-g)x(1-g) \in
I, \quad y_1 = (1-g)y(1-g) \in I_0.$$
Then $p = b_1 + x_1$ and $\|x_1-y_1\| \le \|x-y\| < \ep/2$.
Now, $p$ and $g$ commute, $g$ and
$py_1^2p$ belong to $pI_0p$, $g$ is full in $I_0$, and
$py_1^2p \perp g$. By pure infiniteness of $A$ we deduce that $py_1^2p
\precsim g$. We can now use Lemma~\ref{lm:apprunit}
to conclude that there is a projection $r \in pI_0 p \subseteq pIp$ such that
$\|y_1(p-r)\| < \ep/2$. Hence
$\|x_1(p-r)\| < \ep$.

Now,
\begin{eqnarray*}
p-r & =  & p^*(p-r)p \\ & = &
b_1(p-r)b_1 + b_1(p-r)x_1 + x_1(p-r)b_1 + x_1(p-r)x_1 \\ & = & b_2 + x_2,
\end{eqnarray*}
where
$$b_2 = b_1(p-r)b_1 \in B, \qquad x_2 = b_1(p-r)x_1 + x_1(p-r)b_1 +
x_1(p-r)x_1 \in I.$$
Note that
\begin{eqnarray*}
\|x_2\| & \le & \|b_1\|\|(p-r)x_1\| + \|x_1(p-r)\|\|b_1\| + \|x_1(p-r)\|^2
\\ & \le & 2\|b\|\|x_1(p-r)\| + \|x_1(p-r)\|^2 \\  & \le &
2\|b\|\ep + \ep^2 \; < \; 1/2,
\end{eqnarray*}
where it has been used that $x_1$ is self-adjoint. This shows that the distance
from $b_2$ to the projection $p-r$ is less than $1/2$, whence $1/2$ is not
in the spectrum of $b_2$. The function $f = 1_{[1/2,\infty)}$ restricts to a
continuous function on $\spek(p-r)$ and on $\spek(b_2)$, whence
$$p-r = f(p-r) = f(b_2) + x_3$$
for some $x_3 \in I$. We can take $q$ to be $f(b_2)$.
\end{proof}

\begin{lemma} \label{lm:K0}
Let $A$ be a separable \pif{} \Cs{} with property (IP). Then
$$K_0(A) = \{ [p] : p \; \; \text{is a projection in} \; \; A\}.$$
\end{lemma}

\begin{proof} By Proposition~\ref{prop:stable} every element in
  $K_0(A)$ is represented by a difference $[p_0]-[q_0]$, where
  $p_0,q_0$ are projections in $A \otimes \cK$. Upon replacing $p_0$
  and $q_0$ with $p_0 \oplus q_0$ and $q_0 \oplus q_0$, respectively,
  we can assume that $q_0$ belongs to the ideal
  generated by $p_0$, whence $q_0 \sim q_1 \le p_0$ for some
  projection $q_1$ by pure infiniteness
  of $A$. The projection $p_0-q_1 \in A \otimes \cK$ is equivalent to
  a projection $p \in A$ by Lemma~\ref{lm:her}; and $[p_0]-[q_0] =
  [p_0]-[q_1] = [p_0-q_1] = [p]$.
\end{proof}

\begin{lemma} \label{lm:lift0}
Let
$$\xymatrix{0 \ar[r] & I \ar[r] & A \ar[r]^-\pi & B \ar[r] & 0}$$
be an extension where $A$ is a separable \pif{} \Cs{} with property
(IP). Let $q$ be a projection in $B$ such that $[q]$ belongs to
$K_0(\pi)(K_0(A))$. Then $A$ contains an ideal
$A_0$, which is generated by a single projection, such that $q \in
\pi(A_0)$ and $[q] \in K_0(\pi|_{A_0})(K_0(A_0))$.
\end{lemma}

\begin{proof} By Corollary~\ref{cor:ip} there is an increasing net
  $\{A_\alpha\}_\alpha$ of ideals in $A$, each generated
  by a single projection, such
  that $\bigcup_\alpha A_\alpha$ is dense in $A$. By the assumption
  that $[q] \in K_0(\pi)(K_0(A))$,
  and by Lemma~\ref{lm:K0}, there is a projection $r \in A$ such that
  $[\pi(r)] = [q]$. Now, $r \in A_{\alpha_1}$ and $q \in
  \pi(A_{\alpha_2})$ for suitable $\alpha_1$ and $\alpha_2$. We can
  therefore take $A_0$ to be $A_\alpha$, when $\alpha$ is chosen greater
  than or equal to both $\alpha_1$ and $\alpha_2$.
\end{proof}

\noindent In the lemma below we identify $A$ with the upper left
corner of $M_n(A)$, and thus view $A$ as a hereditary sub-\Cs{} of
$M_n(A)$ for any $n \in \N$.

\begin{lemma} \label{lm:lift1} Let
$$\xymatrix{0 \ar[r] & I \ar[r] & A \ar[r]^-\pi & B \ar[r] & 0}$$
be an extension where $A$ is a separable \pif{} \Cs{} with property
(IP). Then a full projection $q$ in $B$ lifts to a projection in $A
+ M_4(I)$ if and only if $[q] \in K_0(\pi)(K_0(A))$.
\end{lemma}

\begin{proof} The pre-image of $B \subseteq M_4(B)$ under the quotient
  mapping $\pi \otimes \Id_{M_4} \colon M_4(A) \to M_4(B)$ is $A +
  M_4(I)$. Hence it suffices to show that $q$
  lifts to a projection $p \in M_4(A)$.

By Lemmas~\ref{lm:K0} and \ref{lm:lift0}, possibly upon replacing
$A$ by an ideal in $A$, we can assume that $A$
contains a full projection $e$ and a (not necessarily full)
projection $p_1$ such that $[\pi(p_1)] = [q]$ in $K_0(B)$. Since $e$ is full and 
properly infinite there are mutually orthogonal subprojections $e_0$ and
$e_1$ of $e$ such that $e_0$ is full in $A$, $[e_0] = 0$ in
$K_0(A)$, and $e_1 \sim p_1$. Set $p' = e_0+e_1$. Then
  $[\pi(p')] = [q]$ in $K_0(B)$, and $\pi(p')$ and $q$ are both full
  and properly infinite in $B$, so they are equivalent (by
  \cite[Theorem~1.4]{Cuntz:KOn}). It follows that $\pi(p')$ is homotopic to
  $q$ inside $M_4(B)$; and by standard non-stable $K$-theory, see e.g.\
\cite[Lemma 2.1.7, Proposition 2.2.6, and 1.1.6]{RorLarLau:k-theory}, we
conclude that $q$ lifts to a projection $p$ in $M_4(A)$.
\end{proof}

\noindent Using pure infiniteness of $A$ one can improve the lemma
above to get the lifted projection inside $A + M_2(I)$ (instead of in $A +
M_4(I)$). However, one cannot always get the lift in $A + I$ as
Example~\ref{ex:lift} below shows. First we state and prove our main
lifting result for projections in \pif{} \Cs s with the ideal property:

\begin{proposition} \label{prop:lift}
Every separable, \pif, $K_0$-liftable \Cs{} $A$ with property
(IP) has the following projection lifting property: For any hereditary
sub-\Cs{} $A_0$ of $A$ and for any ideal $I_0$ in
$A_0$, every projection in the quotient $A_0/I_0$ lifts to a
projection in $A_0$.
\end{proposition}

\begin{proof} Let $A_0$ and $I_0$ be as above, and let $q$ be a projection in
  $A_0/I_0$. We must show that $q$ lifts to a projection in
  $A_0$. Let $\pi \colon A_0 \to A_0/I_0$ denote the quotient mapping.
  Upon passing to a hereditary sub-\Cs{} of $A_0$ (the
  pre-image $\pi^{-1}(q(A_0/I_0)q)$) we can assume that $q$ is full in
  $A_0/I_0$ (and even that $q$ is the
  unit for $A_0/I_0$). By Lemma~\ref{lm:lift0} (and
  Proposition~\ref{prop:ipher}), possibly upon replacing
  $A_0$ with an ideal of $A_0$, we can further assume that
  $A_0$ contains a full projection, say $g$ (and that $q \in
  \pi(A_0)$).

Put $A_{00} = (1-g)A_0(1-g)$ and $I_{00} = A_{00} \cap I_0$.  It
follows from Lemma~\ref{lm:K0onto} and our assumption that the map
$K_0(A_{00}) \to K_0(A_{00}/I_{00})$ is onto.
We can now use Lemma~\ref{lm:lift1} to lift $q-\pi(g)$ to a
projection $p'$
in $A_{00}+M_4(I_{00})$. Thus $p'' = p'+g \in A_0 + M_4(I_0)$ is a lift
of $q$, and $gA_0g \subseteq p''M_4(A_0)p'' \cap A_0$ is full in $A_0$.
As $A$ is assumed to have property (IP) we obtain from
Proposition~\ref{prop:ipher} that $I_0$ has property (IP), and so we
can use Lemma~\ref{lm:BP} to get a projection $p \in A_0$ such that $p-p''
\in M_4(I_0)$; and $p$ is a lift of $q$.
\end{proof}

\begin{example} \label{ex:lift}
Consider the \Cs{}
$$A = \{ f \in C([0,1],\cO_2) : f(1) = sf(0)s^*\},$$
where $s \in \cO_2$ is any non-unitary isometry. Let $\pi \colon A \to
\cO_2$ be given by $\pi(f) = f(0)$. Then we have a short exact sequence
$$\xymatrix{0 \ar[r] & C_0((0,1),\cO_2) \ar[r] & A \ar[r]^-\pi & \cO_2
  \ar[r] & 0.}$$
The map $K_0(A) \to K_0(\cO_2)$ is surjective, because
$K_0(\cO_2)=0$. (One can show that $A \cong A \otimes \cO_2$, and
hence that $A$ is $K_0$-liftable, cf.\ Lemma~\ref{lm:K0lift-tensor}.) 
However, the unit $1 \in \cO_2$
does not lift to a projection
in $A$, because $1$ is not homotopic to $ss^* \ne 1$ inside $\cO_2$.

Of course, the ideal $C_0((0,1),\cO_2)$ does not have property (IP), so
this example does not contradict
Proposition~\ref{prop:lift}. But the example does show that
Proposition~\ref{prop:lift} is false without
the assumption that $A$ (and hence the ideal $I_{00}$)
has property (IP), and it shows that
Lemma~\ref{lm:lift1} does not hold with $A+M_4(I)$ replaced with $A+I$.
\end{example}

\section{The main result} \label{sec:main}

\noindent Here we state and prove our main result described in the
abstract. Let us set up some notation.

Let $g_\ep \colon \R^+ \to \R^+$ be the continuous function
$$g_\ep(t) = \begin{cases} (\ep-t)/\ep, & t \le \ep\\0, & t \ge
  \ep \end{cases}.$$
If $A$ is a non-unital \Cs{} and $a$ is a positive
element in $A$, then $g_\ep(a)$ belongs to the unitization of $A$, but not
to $A$. However,
$$I_\ep (a) := \overline{Ag_\ep (a)A}$$
is an ideal in $A$, and
$$H_\ep (a) := \overline{g_\ep(a)Ag_\ep(a)}$$
is a hereditary subalgebra of $A$. The hereditary sub-\Cs{} $H_\ep(a)$ is
full in $I_\ep(a)$, i.e., $I_\ep(a) = \overline{A H_\ep(a) A}$.

The quotient \Cs{} $A/I_\ep(A)$ is unital and $a + I_\ep(a)$ is
invertible in $A/I_\ep(a)$ (provided that $I_\ep(a)$ is different from
$A$). Indeed, $h(a) + I_\ep(a)$ is a unit for $A/I_\ep(a)$ and $f(a) +
I_\ep(a)$ is the inverse to $a + I_\ep(a)$, when
$$h(t) = \begin{cases} \ep^{-1}t, & t \le \ep\\ 1, & t \ge \ep
\end{cases},
\qquad f(t) = \begin{cases} \ep^{-2}t, & t \le \ep \\ 1/t, & t \ge
  \ep\end{cases}.$$

\begin{lemma} \label{lm:2}
Assume that $A$ is a separable \pif{} \Cs{} whose primitive ideal
space has a basis of compact-open sets. Let $a$ be a positive
element in $A$ and let $\ep > 0$. Assume that $I_\ep(a) \ne A$. Then
there is a projection $e$ in $H_\ep(a)$ and an ideal
$I$ in $A$ such that $I = \overline{AeA} \subseteq I_\ep(a)$, and
$A/I$ contains a projection $f$ which is a unit for the element
$(a-\ep)_+ + I$ in $A/I$.
\end{lemma}

\begin{proof} If $I_\ep(a)$ itself were compact, i.e., generated by a single
  projection, then, by Lemma~\ref{lm:her}, it would be generated by a
  projection $e \in H_\ep(a)$. We could then take $I$ to be $I_\ep(a)$ and
  the projection $f$ to be the unit of $A/I_\ep(a)$.

Let us now consider the general case, where $I_\ep(a)$ need not be
compact. Find an increasing net of ideals $I_\alpha$
in $A$, each of which is generated by a single projection, such
that $\bigcup_\alpha I_\alpha$ is dense in $I_\ep(a)$, cf.\
Corollary~\ref{cor:ip}.  Then, for each
$\alpha$, we have a commutative diagram:
$$\xymatrix@C+1pc@R+1pc{ & A/I_\alpha \ar[d]^-{\nu_\alpha} \\ A
  \ar[ur]^-{\pi_\alpha} \ar[r]_-\pi & A/I_\ep(a)}$$
and $\|\pi_\alpha(x)\| \to \|\pi(x)\|$ for all $x \in A$. We saw above that
$\pi(h(a))$ is a unit for $A/I_\ep(a)$; so
$$\lim_\alpha \|\pi_\alpha\big(h(a)-h(a)^2\big)\| = \|\pi\big(h(a) -
h(a)^2\big)\| = 0.$$
We can therefore take $\alpha$ such that
$\|\pi_\alpha\big(h(a)-h(a)^2\big)\| < 1/4$, in which case $1/2$ does not
belong to the spectrum of $\pi_\alpha(h(a))$.

The ideal $I_\alpha$ is by assumption generated by a projection, say $g$;
and as $g$ belongs to $I_\ep(a)$ it is equivalent to a projection $e \in
H_\ep(a)$ by Lemma~\ref{lm:her}; whence $I:= I_\alpha$ is generated by $e$.

The characteristic function $1_{[1/2,\infty)}$ is continuous on the
spectrum of $\pi_\alpha(h(a))$; and it extends to a continuous
function $\varphi \colon \R^+ \to [0,1]$ which satisfies
$\varphi(0)=0$ and $\varphi(1)=1$. Put
$$f = 1_{[1/2,\infty)}\big(\pi_\alpha(h(a))\big) =
\pi_\alpha\big((\varphi \circ h)(a)\big) \in A/I.$$
Then $f$ is a projection, and as $(\varphi \circ h)(a) \gange (a-\ep)_+ =
(a-\ep)_+$, we have
$$f \gange \pi_\alpha\big((a-\ep)_+\big) =
\pi_\alpha\big((\varphi \circ h)(a)\gange (a-\ep)_+\big) =
\pi_\alpha\big((a-\ep)_+\big),$$
as desired.
\end{proof}

\begin{theorem} \label{thm:rr0}
Let A be a separable \pif{} \Cs{}. Then the real rank of $A$ is zero
if and only if $A$ is $K_0$-liftable (cf.\ Definition~\ref{def:K0onto}) and
the primitive ideal space of $A$ has a basis for its topology consisting of
compact-open sets. 
\end{theorem}

\begin{proof} If $\RR(A)=0$, then $A$ has property (IP), whence
  $\Prim(A)$ has a basis consisting of compact-open sets, cf.\
  Proposition~\ref{prop:IP}.  As
  remarked below Definition~\ref{def:K0onto}, it follows from
  \cite{BroPed:realrank} that every \Cs{} of real rank
  zero is $K_0$-liftable. This proves the ``only if'' part.

We proceed to prove
  the ``if'' part, and so we assume that $A$ is $K_0$-liftable and
  that $\Prim(A)$ has a basis of compact-open sets. Then, by
  Proposition~\ref{prop:IP}, $A$ has property (IP).

To show that $\RR(A)=0$ we show that each hereditary sub-\Cs{} of
$A$ has an approximate unit consisting of projections. Hereditary
sub-\Cs s of \pif{} \Cs s are again \pif{} (see \cite{KirRor:pi}),
and it follows from Lemma~\ref{lm:K0onto} and
Proposition~\ref{prop:ipher} that any hereditary sub-\Cs{} of $A$ is
$K_0$-liftable and has property (IP). Upon replacing a hereditary
sub-\Cs{} of $A$ by $A$ itself, it suffices to show that $A$ has an
approximate unit consisting of projections. To this end it suffices
to show that, given a positive element $a$ in $A$ and $\ep > 0$,
there is a projection $p$ in $A$ such that $\|a-ap \| \le 3\ep$.

Let $I_\ep(a)$ and $H_\ep(a)$ be as defined above Lemma~\ref{lm:2}. Then,
as already observed, $H_\ep(a)$ is a full hereditary sub-\Cs{} of
$I_\ep(a)$; and $\|ae\| \le \ep
\|e\|$ for all $e \in H_\ep(a)$ by construction of $H_\ep(a)$. 

Suppose that $I_\ep(a) = A$. Then $H_\ep(a)$ is a full hereditary sub-\Cs{}
in $A$.
By Lemma~\ref{lm:proj2} there is a
projection $f$ in $A$ with $(a-\ep)_+ \precsim f$; and by
Lemma~\ref{lm:her}, $f$ is equivalent to a projection $e \in
H_\ep(a)$. As $(a-\ep)_+ \precsim e$ and $(a-\ep)_+ \perp
e$ we can use Lemma~\ref{lm:apprunit} to find a projection $p
\in A$ such that $\|(a-\ep)_+(1-p)\| \le \ep$, whence $\|a(1-p)\| \le
2 \ep \le 3 \ep$.

Suppose now that $I_\ep(a) \ne A$. Let $e \in H_\ep(a)$, $I =
\overline{AeA}$, and $f \in A/I$ be as in Lemma~\ref{lm:2}, and let
$\pi \colon  A \to A/I$ denote the quotient mapping. Note that
$\pi\big((1-e)A(1-e)\big) = \pi(A)$. It follows from
Proposition~\ref{prop:lift} that $f$ lifts to a projection $q$ in
$(1-e)A(1-e)$. Consider the element $x= (a-\ep)_+(1-e-q) =
(a-\ep)_+(1-q)$, which belongs to  
$I$ because $\pi(x)=0$. Hence $x^*x \precsim e$ by pure infiniteness
of $A$, and $x^*x$ is clearly
orthogonal to $e$. As both $e$ and $x^*x$ belong to the corner \Cs{} 
$(1-q)A(1-q)$ and the relation $x^*x \precsim e$ also holds relatively
to this corner, it follows from Lemma~\ref{lm:apprunit} 
that there is a projection $r$ in $(1-q)A(1-q)$ such that
$\|x^*x(1-r)\| \le \ep^2$, whence 
$$\|(a-\ep)_+(1-e-q)(1-r)\| = \|x(1-r)\| \le \|x^*x(1-r)\|^{1/2} \le \ep.$$
Recall that $e \perp q$ and $r \perp q$. Put $p = r+q$, and note
that $(1-e-q)(1-r) = (1-e)(1-q)(1-r) = (1-e)(1-p)$. We can now deduce that
\begin{eqnarray*}
\|a(1-p)\| & \le & \|a(1-e)(1-p)\| + \|ae(1-p)\| \\
& \le & \|a(1-e-q)(1-r)\| + \|ae\| \\
& \le & 2\ep + \ep  \; = \; 3\ep.
\end{eqnarray*}
\end{proof}

\noindent Our theorem above generalizes, in the separable case,
Zhang's theorem (from \cite{Zhang:infsimp}) that all simple, \pif{}
\Cs s are of real rank zero. The primitive ideal space of a simple
\Cs{} consists of one point (the $0$-ideal) and hence trivially has
a basis of compact-open sets, and any simple \Cs{} is automatically
$K_0$-liftable (as remarked below Definition~\ref{def:K0onto}).

\newpage

\begin{corollary} \label{cor:O2}
Let $A$ be any separable \Cs.
\begin{enumerate}
\item $\RR(A \otimes \cO_\infty) = 0$ if and only if $\Prim(A)$ has a
  basis consisting of compact-open sets and $A$ is $K_0$-liftable.
\item  The following three conditions are equivalent:
\begin{enumerate}
\item $\RR(A \otimes \cO_2) = 0$,
\item $A \otimes \cO_2$ has property (IP),
\item $\Prim(A)$ has a basis consisting of compact-open sets.
\end{enumerate}
If, in addition, $A$ is \pif, then conditions {\rm{(a)--(c)}} above are
equivalent to:
\begin{itemize}
\item[{\rm{(d)}}] $A$ has property (IP).
\end{itemize}
\end{enumerate}
\end{corollary}

\begin{proof} The \Cs s $A \otimes \cO_\infty$ and $A \otimes \cO_2$
  are \pif{} and separable (cf.\ \cite{KirRor:pi} and
  Remark~\ref{rm:pi}). The ideal lattices $\Ideal(A)$, $\Ideal(A
  \otimes \cO_2)$, and $\Ideal(A \otimes \cO_\infty)$ are isomorphic,
  cf.\ Lemma~\ref{lm:K0lift-tensor} and its proof, whence---by
  separability---$\Prim(A)$, $\Prim(A \otimes
  \cO_\infty)$ and $\Prim(A \otimes \cO_2)$ are homeomorphic. It
  follows from Lemma~\ref{lm:K0lift-tensor} that $A \otimes \cO_2$ is
  $K_0$-liftable, and that $A \otimes \cO_\infty$ is $K_0$-liftable if and only
  if $A$ is $K_0$-liftable. The claims of the corollary now follow from
  Theorem~\ref{thm:rr0} and Proposition~\ref{prop:IP}.
\end{proof}

\noindent Extensions of separable \Cs s with property (IP) need not
have property (IP) (not even after being tensored by the compacts),
cf.\ \cite{Pas:ideal_prop_AH}. But in the purely infinite case we
have the following:


\begin{corollary} \label{cor:extension}
Let  $0 \to I \to A \to B \to 0$ be an extension of separable \Cs s.
\begin{enumerate}
\item If $\Prim(I)$ and $\Prim(B)$ have basis for their topology consisting
  of compact-open sets, then so does $\Prim(A)$. 
\item If $I$ and $B$ are \pif{} and with property (IP), then so is $A$.
\end{enumerate}
\end{corollary}

\begin{proof} (i).  It suffices to show that $\RR(A \otimes \cO_2) = 0$, cf.\
  Corollary~\ref{cor:O2}. But
$$\xymatrix{0 \ar[r] & I \otimes \cO_2 \ar[r] & A \otimes \cO_2 \ar[r]
  & B \otimes \cO_2 \ar[r] & 0}$$
is exact, because $\cO_2$ is exact, $\RR(I \otimes \cO_2) = 0$, $\RR(B
\otimes \cO_2) = 0$ by  Corollary~\ref{cor:O2}, and $K_1(I \otimes
\cO_2) = 0$. It therefore follows from \cite[Theorem~3.14 and
Proposition~3.15]{BroPed:realrank} that $\RR(A \otimes \cO_2)$ is zero.

(ii). It follows from (i) and Corollary~\ref{cor:O2} that $\RR(A \otimes
\cO_2) = 0$, whence $A$ has property (IP), again by
Corollary~\ref{cor:O2}. It is shown in \cite[Theorem~4.19]{KirRor:pi} that
extensions of \pif{} \Cs s again are \pif. 
\end{proof}

\noindent It is shown in \cite[Proposition~2.16]{BlanKir:pi3} that
$\Prim(A \otimes B)$ is homeomorphic to $\Prim(A) \times \Prim(B)$
when either $A$ or $B$ is exact. It follows in particular that
$\Prim(A \otimes B)$ has a basis of compact-open sets if both
$\Prim(A)$ and $\Prim(B)$ have basis of compact-open sets and if and
one of $A$ and $B$ is exact.   

The tensor product $A \otimes B$ can contain unexpected ideals if both
$A$ and $B$ are non-exact. More specifically, it
follows from a theorem of Kirchberg that if $C$ is a simple \Cs{} and
$H$ is an infinite-dimensional (separable) Hilbert space, then $B(H)
\otimes C$ has more than the three obvious ideals (counting the two
trivial ones) 
if and only if $C$ is non-exact. Part (i) of
the proposition below shows that $\Prim(A \otimes B)$ can be much
larger than $\Prim(A) \times \Prim(B)$. 

\begin{proposition} \label{prop:tensor-nonexact}
There are separable (necessarily non-exact) \Cs s $A$ and $C$ such that
$\Prim(A)$ consists of 
two points (i.e., $A$ is an extension of two simple \Cs s) and $\Prim(C)$
consists of one point (i.e., $C$ is simple) such that:
\begin{enumerate}
\item $\Prim(A \otimes C)$ does not have a basis for its topology
  consisting of compact-open sets; in particular, $\Prim(A \otimes C)$ is
  infinite. 
\item The \Cs s $A \otimes \cO_2$ and $C \otimes \cO_2$ are \pif{} and
  of real rank zero (and hence with property (IP)), 
  but their tensor product $(A \otimes \cO_2) \otimes (C \otimes \cO_2)$
  does not have property (IP) (and hence is not of real rank zero). 
\end{enumerate}
\end{proposition}

\begin{proof} Let $C$ be the non-exact, simple, unital, separable \Cs{} with
  stable rank one and real rank zero constructed by Dadarlat in
  \cite{Dad:non-exact} (see also \cite[2.1]{PasRor:tensIP}). Let $A$ be the
  (also non-exact) separable sub-\Cs{} of $B(H)$ constructed 
  in \cite[Theorem~2.6]{PasRor:tensIP}. Then $A \otimes C$, and hence also
  $A \otimes C \otimes \cO_2$, contain more than three ideals (including the
  two trivial ones) (by \cite[Theorem~2.6]{PasRor:tensIP}). 

It follows from \cite[Proposition~2.2]{PasRor:tensIP} (following Dadarlat's
construction) that there is a UHF-algebra $B$ which is shape
equivalent to $C$, whence the following holds: For any \Cs{} $D$, the
subsets of $\Ideal(B \otimes D)$ and of $\Ideal(C \otimes D)$,
consisting of all ideals that are generated by projections, are order
isomorphic.   

The ideal lattice of $A \otimes B \otimes \cO_2$ is order isomorphic to the
ideal lattice of $A$ (because $B \otimes \cO_2$ is simple and exact), so $A \otimes B
\otimes \cO_2$ has three ideals (including the two trivial ideals),
and each of these three ideals is generated by its projections. It
follows that $A \otimes C \otimes \cO_2$ also has precisely three
ideals that are generated by 
projections. Hence $A \otimes C \otimes \cO_2$ has at least one ideal which
is not generated by projections. We conclude that $A \otimes C \otimes
\cO_2$ does not have property (IP). Hence $\Prim(A \otimes C)$ does not
have a basis of compact-open sets (by Corollary~\ref{cor:O2}) and $(A
\otimes \cO_2) \otimes (C \otimes \cO_2)$, which is isomorphic to $A
\otimes C \otimes \cO_2$, does not have property (IP). It follows from
Corollary~\ref{cor:O2} that $A \otimes \cO_2$ and $C \otimes \cO_2$ both
are of real rank zero. 
\end{proof}

\begin{proposition} \label{prop:IPtensor}
    Let $A$ and $B$ be \Cs{}s with property (IP).  Assume that $A$ is
    exact and that $B$ is purely infinite.  Then $A\otimes B$ is
    purely infinite and with property (IP).
\end{proposition}

\begin{proof}
    Since $B$ is purely infinite and with property (IP),
    Proposition~\ref{prop:3pi} implies that $B$ is strongly purely
    infinite.  But a recent result of Kirchberg says that
    if $C$ and $D$ are \Cs{}s such that one of $C$ or $D$ is exact and the
    other is strongly purely infinite, then $C\otimes D$ is
    strongly purely infinite (see \cite{Kir:Oberwol}).  Hence, by this
    result of Kirchberg it follows that $A\otimes B$ is strongly
    purely infinite, and hence purely infinite.  Also, since $A$ is
    exact and $A$ and $B$ have property (IP), by
    \cite[Corollary~1.3]{PasRor:tensIP} (based on another result of
    Kirchberg), it follows that $A\otimes B$ has property (IP).
\end{proof}

\noindent There are well-known examples of two separable nuclear \Cs s each
of real rank zero whose minimal tensor product is a \Cs{} not of
real rank zero (see \cite{KodOsa:RR0_tensor2}). This phenomenon is 
eliminated when tensoring with $\cO_2$: 

\begin{corollary} \label{cor:tensor2}
    Let $A$ and $B$ be separable \Cs s with property (IP) (or of real
    rank zero).  Assume that $A$ is exact.  Then $A\otimes B \otimes \cO_2$
    is of real rank zero.
\end{corollary}

\begin{proof} It follows from Proposition~\ref{prop:IPtensor} that
  $A\otimes B \otimes \cO_2$ is \pif{} and with property (IP), whence
  this \Cs{} is of real rank zero by Corollary~\ref{cor:O2}.
\end{proof}

\noindent The two conditions (on the primitive ideal space and on
$K_0$-liftability)  in Theorem~\ref{thm:rr0} are
independent. There are \pif{} \Cs s that are $K_0$-liftable while others
are not, and there are \pif{} \Cs s whose primitive ideal
space has a basis of compact-open sets, and others where this does not
hold. All four combinations exist. The \Cs s $C([0,1]) \otimes
\cO_\infty$ and $C([0,1]) \otimes \cO_2$ are purely infinite with
primitive ideal space homeomorphic to $[0,1]$, and this space does not have a
basis of compact-open sets (i.e., is not totally disconnected); the
latter \Cs{} is $K_0$-liftable and the former is not (consider the
surjection $C([0,1]) \otimes \cO_\infty \to C(\{0,1\}) \otimes
\cO_\infty$). More examples are given below:

\begin{example}[The case where the primitive ideal space is finite]
  \label{ex:finite} Every subset of a finite T$_0$-space is compact
  (has the Heine-Borel property), so if $A$ is a \Cs{} for which
  $\Prim(A)$ is finite, then $\Prim(A)$ has a basis of compact open
  sets. Suppose that $\Prim(A)$ is finite and that $A$ is \pif. Then
  $\Ideal(A)$ is a
  finite lattice, and there exists a decomposition series
$$0 = I_0 \lhd I_1 \lhd I_2 \lhd \cdots \lhd I_n = A,$$
where each $I_j$ is a closed two-sided ideal in $A$, and where
each successive quotient $I_j/I_{j-1}$, $j=1,2, \dots, n$, is simple.

It follows from Theorem~\ref{thm:rr0} that $A$
is of real rank zero if and only if $A$ is $K_0$-liftable (when $A$ is
separable). Actually, one can obtain this result (also in the
non-separable case) from Zhang's theorem,
which tells us that $I_j/I_{j-1}$ is of real rank zero for all $j$,
being simple and \pif, and from Brown and Pedersen's extension
result in \cite[Theorem~3.14 and Proposition~3.15]{BroPed:realrank},
applied to the extension
$$\xymatrix{0 \ar[r] & I_{j-1} \ar[r] & I_j \ar[r]
  & I_j/I_{j-1} \ar[r] & 0,}$$
which yields that $\RR(I_j)=0$ if (and only if) $\RR(I_{j-1})=0$ and
$K_0(I_j) \to K_0(I_j/I_{j-1})$ is surjective. Hence $\RR(A) = 0$ if and
only if $K_0(I_j) \to K_0(I_j/I_{j-1})$ is surjective for all $j=1,2, \dots,
n$. The latter is equivalent to $A$ being $K_0$-liftable (as one
easily can deduce from Lemma~\ref{lm:K0lift-extension}). 

In the case where $n=2$ we have an extension $0 \to I \to A \to B
\to 0$, where $I$ and $B$ are \pif{} \Cs s. Here $\RR(A) = 0$ if and
only if the map $K_0(A) \to K_0(B)$ is surjective, or equivalently,
if and only if the index map $\delta \colon K_0(B) \to K_1(I)$ is
zero. Let $G_0, G_1, H_0, H_1$ be arbitrary countable abelian groups
and let $\delta \colon G_0 \to H_1$ be any group homomorphism. Then
there are stable Kirchberg algebras $I$ and $B$ in the UCT-class
such that $K_j(B) \cong G_j$ and $K_j(I) \cong H_j$, and an
essential extension $0 \to I \to A \to B \to 0$ whose index map
$K_0(B) \to K_1(I)$ is conjugate to $\delta$. 

In particular, if $G_0$, $H_1$, and $\delta$ are chosen such that
$\delta$ is non-zero, then $A$ is not $K_0$-liftable and hence not of real
rank zero; but $A$ is $K_0$-liftable and of real rank zero whenever
$\delta$ is zero. Evidently, both situations can occur.

Let us finally note that $\Prim(A)$, if finite,
is Hausdorff if and only if the
topology on $\Prim(A)$ is the discrete topology, which happens if and only
if $A$
is the direct sum of $n$ simple \pif{} \Cs s. Here, $K_0$-liftability
is automatic. Note also that $\Prim(A)$ is totally disconnected
(meaning that all connected components are singletons) if and only if
$\Prim(A)$ is Hausdorff.
\end{example}

\noindent\textbf{Acknowledgments:} The first named author was
partially supported by NSF Grant DMS-0101060 and also by a FIPI grant
from the University of Puerto Rico.

\providecommand{\bysame}{\leavevmode\hbox to3em{\hrulefill}\thinspace}
\providecommand{\MR}{\relax\ifhmode\unskip\space\fi MR }
\providecommand{\MRhref}[2]{%
  \href{http://www.ams.org/mathscinet-getitem?mr=#1}{#2}
}
\providecommand{\href}[2]{#2}

\vspace{.3cm}
\noindent{\sc Department of Mathematics, University of Puerto Rico,
  P.O.\ Box 23355, San Juan, Puerto Rico 00931, USA}\\

\vspace{.2cm}
\noindent{\sl E-mail address:} {\tt cpasnic@upracd.upr.clu.edu}\\

\vspace{.5cm}

\noindent{\sc Department of Mathematics, University of Southern
  Denmark, Odense,
  Campusvej~55, 5230 Odense M, Denmark}

\vspace{.3cm}

\noindent{\sl E-mail address:} {\tt mikael@imada.sdu.dk}\\
\noindent{\sl Internet home page:}
{\tt www.imada.sdu.dk/$\,\widetilde{\;}$mikael/welcome} \\

\end{document}